\documentclass{amsart}
\usepackage{amssymb}
\begin{document}
\newtheorem{theorem}{Theorem}[section]
\newtheorem{lemma}[theorem]{Lemma}
\newtheorem{definition}[theorem]{Definition}
\def\operatorname#1{{\rm#1\,}}
\def\qedbox{\hbox{$\rlap{$\sqcap$}\sqcup$}}
\def\range{\operatorname{range}}
\def\Pspan{\operatorname{span}}
\def\RR{\mathcal{R}}
\def\rank{\operatorname{rank}}
\def\id{\operatorname{Id}}
\def\trace{\operatorname{trace}}
\def\imag{\operatorname{Im}}
\font\pbglie=eufm10
\def\SS{{\text{\pbglie S}}}
\def\spec{\operatorname{Spec}}
 \makeatletter
  \renewcommand{\theequation}{%
   \thesection.\arabic{equation}}
  \@addtoreset{equation}{section}
 \makeatother
\def\pix{\rho_{\mathcal{X}}}
\title[Higher order Jordan Osserman manifolds]{Higher order Jordan Osserman\\ Pseudo-Riemannian manifolds}
\author{Peter B. Gilkey, Raina
Ivanova, and Tan Zhang}
\begin{address}{PG: Mathematics Department, University of Oregon, Eugene OR 97403 USA.\newline\phantom{...........} email:
gilkey@darkwing.uoregon.edu}\end{address}
\begin{address}{RI: Mathmatics Department,
University of Hawaii - Hilo,
200 W. Kawili St.,\newline\phantom{..........}
Hilo  HI 96720 USA. email: rivanova@hawaii.edu}\end{address}
\begin{address}{TZ: Department of Mathematics and Statistics,
Murray State University, Murray\newline\phantom{...........} KY 42071 USA
email: tan.zhang@murraystate.edu}\end{address}
%\date{19 May 2002 Version 1-z}
\begin{abstract} We study the higher order Jacobi
operator in pseudo-Riemannian geometry. We exhibit a family of manifolds so that this operator has constant
Jordan normal form on the Grassmannian of subspaces of signature $(r,s)$ for certain values of $(r,s)$. These
pseudo-Riemannian manifolds are new and non-trivial examples of higher order Osserman manifolds.
\end{abstract}
\maketitle
\font\pbglie=eufm10
\def\Gr{\text{Gr}}

\section{Introduction}\label{Sect1}

Let $(M,g)$ be a pseudo-Riemannian manifold of signature $(p,q)$ and dimension $m=p+q$. Let $R(\cdot,\cdot)$
be the associated Riemann curvature operator. The {\it Jacobi operator} $J(X):Y\rightarrow R(Y,X)X$ is a
self-adjoint operator which plays an important role in the study of geodesic variations and in many other applications. 

We say that $(M,g)$ is {\it Riemannian} if $p=0$ and {\it Lorentzian} if $p=1$.
Osserman \cite{refOss} observed that if $(M,g)$ is a
local $2$ point homogeneous Riemannian manifold, then the eigenvalues of the Jacobi operator are constant on the unit sphere bundle of
$M$. He wondered if the converse held, i.e. if the eigenvalues of the Jacobi operator are constant on the unit sphere bundle, does
this imply that $(M,g)$ is a local $2$ point homogeneous space (or equivalently, that $(M,g)$ is either flat or is a rank $1$ symmetric
space). This has been shown to be true if $m\ne16$ by Chi \cite{refChia} and
Nikolayevsky \cite{refNik}. The case $m=16$ is still open.

It is natural to pose this question in the pseudo-Riemannian setting as well. Let $S^\pm(M,g)$ be the pseudo-sphere bundles of unit
spacelike ($+$) and timelike ($-$) vectors for a pseudo-Riemannian manifold $(M,g)$ of signature $(p,q)$. Then
$(M,g)$ is said to be {\it spacelike Osserman} (resp. {\it timelike Osserman}) if the the eigenvalues of
$J(\cdot)$ are constant on
$S^+(M,g)$ (resp. on $S^-(M,g))$. The notions spacelike Osserman and timelike Osserman are equivalent and if $(M,g)$ is
either of them, then
$(M,g)$ is said to be {\it Osserman}. It is known that any Lorentzian Osserman manifold has constant sectional curvature
\cite{refBokanBlazicGilkey,refGra}. On the other hand, if $p\ge2$ and $q\ge2$, then there exist Osserman pseudo-Riemannian manifolds
of signature
$(p,q)$ which are not locally homogeneous \cite{refBBGZ,refGVV}.

In the higher signature setting, unlike in the Riemannian setting, the eigenvalue structure does not determine the Jordan normal form
of a symmetric linear operator. We say that $(M,g)$ is {\it spacelike Jordan Osserman} (resp. {\it timelike Jordan Osserman}) if
the Jordan normal form (i.e. the conjugacy class) of $J(\cdot)$ is constant on $S^+(M,g)$ (resp. on $S^-(M,g)$). The notions
spacelike Jordan Osserman and timelike Jordan Osserman are distinct. While spacelike Jordan Osserman or timelike Jordan Osserman
implies Osserman, the reverse implication fails in general. There is an extensive literature on this question; see
\cite{refGRKVL,refGiA,refGil,refGiSt} for further details.

This paper focuses on the higher order Jacobi operator, which was first defined by Stanilov and Videv 
\cite{refSV} in the Riemannian setting. We consider it in the pseudo-Riemannian setting. Let
$\Gr_{r,s}(M,g)$ be the Grassmannian bundle of non-degenerate subspaces of signature $(r,s)$. We say that a pair $(r,s)$ is {\it
admissible} if 
$\Gr_{r,s}(M,g)$ is non-empty and does not consist of a single point, i.e.
$$0\le r\le p,\quad 0\le s\le q,\quad\text{ and }\quad1\le r+s\le m-1.$$
Let $(r,s)$ be admissible and let
$\{e_1,...,e_{r+s}\}$ be an orthonormal basis for a subspace $\pi\in\Gr_{r,s}(M,g)$. The {\it higher order Jacobi
operator} is defined \cite{refGSV} in the pseudo-Riemannian setting  by:
$$\textstyle J(\pi):=\sum_{1\le i\le r+s}(e_i,e_i)J(e_i).$$
The operator $J(\pi)$ is independent of the particular orthonormal basis chosen. One says that $(M,g)$ is
{\it Osserman of type
$(r,s)$} if the eigenvalues of
$J(\cdot)$ are constant on $\Gr_{r,s}(M,g)$. Note that $(M,g)$ is Osserman of type $(0,1)$ or $(1,0)$ if and only if $(M,g)$ is
Osserman. Similarly, one says that $(M,g)$ is {\it Jordan Osserman of type $(r,s)$} if the Jordan normal
form of
$J(\cdot)$ is constant on $\Gr_{r,s}(M,g)$; $(M,g)$ is Jordan Osserman of type $(0,1)$ (resp. $(1,0)$) if and only if
$(M,g)$ is spacelike (resp. timelike) Jordan Ossersman.

We shall use the following basic duality result \cite{refGiIv,refGSV}:

\begin{theorem}\label{Thm1.1} Let $(M,g)$ be a pseudo-Riemannian manifold of signature $(p,q)$. Let $(r,s)$ be admissible.
\begin{enumerate}
\item If $(M,g)$ is Osserman of type $(r,s)$, then $(M,g)$ is Osserman of type $(\bar r,\bar s)$ for
every admissible $(\bar r,\bar s)$ with $\bar r+\bar s=r+s$ or $\bar r+\bar s=m-r-s$.
\item If $(M,g)$ is Jordan Osserman of type $(r,s)$, then $(M,g)$ is Jordan Osserman of type
$(p-r,q-s)$.
\end{enumerate}
\end{theorem}

In view of Theorem \ref{Thm1.1}, we say that $(M,g)$ is {\it $k$ Osserman} if $(M,g)$ is Osserman of type $(r,s)$ for any (and hence
for all) admissible $(r,s)$ with $r+s=k$ or $r+s=m-k$. It is known  that in either the Riemannian or the Lorentzian settings,  if
$(M,g)$ is $k$ Osserman for some $k$ with $2\le k\le m-2$, then $(M,g)$ has constant sectional curvature \cite{refGiA,refGiSt}.
Thus, we shall consider the higher signature setting and assume $p\ge2$ and $q\ge2$ henceforth.

Let $M=\mathbb{R}^{2p}$ with the usual coordinates $(x,y):=(x_1,...,x_p,y_1,...,y_p)$. Let
$$\mathcal{X}:=\Pspan_{1\le i\le p}\{\partial_i^x\}
\quad\text{and}\quad
\mathcal{Y}:=\Pspan_{1\le i\le p}\{\partial_i^y\}$$ 
define two distributions of $TM$. The splitting
$TM=\mathcal{X}\oplus\mathcal{Y}$ is, of course, just the usual splitting
$T(\mathbb{R}^{2p})=T(\mathbb{R}^p)\oplus T(\mathbb{R}^p)$. This defines a projection
\begin{equation}\pix :T(\mathbb{R}^{2p})\rightarrow\mathcal{X}.\label{eqn1.1}\end{equation}
\begin{definition} Let $S^2(T\mathbb{R}^p)$ be the bundle of symmetric bilinear forms on $T\mathbb{R}^p$.
Let $\psi\in C^\infty(S^2(T\mathbb{R}^p))$ be a symmetric $2$ tensor with components
$\psi_{ij}:=\psi(\partial_i^x,\partial_j^x)$. Define a pseudo-Riemannian metric of balanced (or neutral) signature $(p,p)$ on
$\mathbb{R}^{2p}$ by setting:
$$g_\psi(x,y)(\partial_i^x,\partial_j^y)=\delta_{ij},\quad 
  g_\psi(x,y)(\partial_i^y,\partial_j^y)=0,\quad 
  g_\psi(x,y)(\partial_i^x,\partial_j^x):=\psi_{ij}(x).$$
\end{definition} 

The
coefficients of  $g_\psi$ depend on $x$ but not on $y$. Furthermore, the distribution $\mathcal{Y}$ is totally isotropic
with respect to $g_\psi$. In Section \ref{Sect2}, we will show that:

\goodbreak\begin{lemma}\label{Lem1.3} Let $\psi_{ij/kl}:=\partial_k^x\partial_l^x\psi_{ij}$. Let $Z_\nu$ be vector fields on
$(M,g_\psi)$. Then:
\begin{enumerate}
\item $R(Z_1,Z_2)Z_3=0$ if $Z_1\in\mathcal{Y}$, if $Z_2\in\mathcal{Y}$, or if $Z_3\in\mathcal{Y}$.
\item $R(\partial_i^x,\partial_j^x)\partial_k^x
   =\textstyle-\frac12\sum_l(\psi_{il/jk}+\psi_{jk/il}-\psi_{ik/jl}-\psi_{jl/ik})\partial_l^y$.
\end{enumerate}\end{lemma}

If $A$ is a self-adjoint linear operator, then we may use the inner product to define an associated symmetric bilinear form
$\mathcal{A}$ by setting
$$\mathcal{A}(Z_1,Z_2):=(AZ_1,Z_2).$$
 Conversely, every symmetric bilinear form arises in this way. A self-adjoint operator $A$ is said to be {\it positive
semi-definite} if and only if the associated quadratic form
$\mathcal{A}$ is positive semi-definite.
The Jacobi operator and the associated bilinear form which are defined by $(M,g_\psi)$ are supported on $\mathcal{X}$. If we set
$X_\nu:=\pix Z_\nu$, then:
\begin{eqnarray*}
&&J(Z_1)Z_2=J(X_1)X_2,\quad\text{ and }\quad
(J(Z_1)Z_2,Z_3)=(J(X_1)X_2,X_3).\end{eqnarray*}
Since $J(X_1)X_1=0$, $\rank(J(Z_1))\le p-1$.

\begin{definition}\label{defn1.4}\ \begin{enumerate}
\item If $\psi\in C^\infty(S^2(T\mathbb{R}^p))$ and if $K$ is a compact subset of
$\mathbb{R}^p$, define the semi-norm
$|\psi|_K=\max_{x\in K;1\le i,j,k,l\le p}|\psi_{ij/kl}|(x)$.
\medbreak \item Let $f\in C^\infty(\mathbb{R}^p)$. Define a symmetric $2$ tensor field  $\psi_f$ on $T\mathbb{R}^p$ by setting
$\psi_{f,ij}=\partial_i^xf\cdot\partial_j^xf$. Let $H(f):=\partial_i^x\partial_j^xf$ be the Hessian.
\medbreak\item Let $\Psi$ be the set of all $\psi\in C^\infty(S^2(T\mathbb{R}^p))$ so that
$J(X)$ is positive semi-definite of rank
$p-1$ for every $0\ne X\in T\mathbb{R}^p$.\end{enumerate}\end{definition}

The metrics $g_\psi$ for $\psi\in\Psi$ will be important in our discussion. We show that this class of metrics is
non-trivial by establishing the following result in Section
\ref{Sect2}:

\begin{lemma}\label{Lem1.5} Let $p\ge2$.
\begin{enumerate}
\item If $f\in C^\infty(\mathbb{R}^p)$ and if $H(f)$ is positive definite, then $\psi_f\in\Psi$.
\item $\Psi$ is conelike, i.e. if $0<a\in\mathbb{R}$ and if $\psi\in\Psi$, then $a\psi\in\Psi$.
\item $\Psi$ is convex, i.e. if $t\in[0,1]$ and if $\psi_i\in\Psi$, then $t\psi_1+(1-t)\psi_2\in\Psi$.
\item Let $K$ be a compact subset of $\mathbb{R}^p$, let $\phi\in C_0^\infty(K)$, and let
$\psi_0\in\Psi$. There exists
$\varepsilon=\varepsilon(\phi,\psi_0)>0$ so that if $\psi_1\in C^\infty( S^2(T\mathbb{R}^p))$ and $|\psi_1|_K<\varepsilon$, then
$\psi_0+\phi\psi_1\in\Psi$.
\end{enumerate}
\end{lemma}

If $(M,g)$ has constant sectional curvature, then $(M,g)$ is Jordan Osserman of type $(r,s)$ for any admissible $(r,s)$. 
There are a number of spacelike Jordan Osserman and timelike Jordan Osserman manifolds which are known
\cite{refBBGZ,refGRKVL,refGIZ}. Higher order Osserman tensors have been constructed
\cite{refGiIv,refGSV,refSt}. Bonome, Castro, and Garcia-Rio have exhibited higher order Osserman manifolds \cite{refBCG} in
signature $(2,2)$. However, no examples of higher order Jordan Osserman manifolds were known prevously which did
not have constant sectional curvature. In Section
\ref{Sect4}, we will demonstrate examples of pseudo-Riemannian manifolds which are Jordan Osserman of type
$(r,s)$ for certain, but not for all values of $(r,s)$:

\begin{theorem}\label{Thm1.6} Let $p\ge2$.\ 
\begin{enumerate}
\item $(M,g_\psi)$ is $k$ Osserman for every admissible $k$.
\item If $\psi\in\Psi$, then $(M,g_\psi)$ is:
\begin{enumerate}
\item Jordan Osserman of types $(0,r)$, $(p,p-r)$, $(r,0)$, $(p-r,p)$ if $0<r\le p$;
\item not Jordan Osserman of type $(r,s)$ if $0<r<p$ and $0<s<p$.
\end{enumerate}\end{enumerate}
\end{theorem}

Theorem \ref{Thm1.6} deals with the balanced signature $p=q$. Let $\mathbb{R}^{(u,v)}$ denote Euclidean space
with the usual inner product of signature $(u,v)$. In Section \ref{Sect5}, we construct examples of
pseudo-Riemannian manifolds of more general signatures by taking the isometric product of $(M,g_\psi)$ with
$\mathbb{R}^{(u,v)}$:

\begin{theorem}\label{Thm1.7}Let $p\ge2$ and let $(N,g_N):=(M,g_\psi)\times\mathbb{R}^{(u,v)}$ 
be the isometric product pseudo-Riemannian manifold of signature $(\bar p,\bar q):=(p+u,p+v)$.
Then:\begin{enumerate}
\item $(N,g_N)$ is $k$ Osserman for every admissible $k$. 
\item If $\psi\in\Psi$, then:
\begin{enumerate} 
\item $(N,g_N)$ is Jordan Osserman
\begin{enumerate}
\item of types $(r,0)$ and $(\bar p-r,\bar q)$ if $u=0$ and if $0<r\le p$;
\item of types $(0,s)$ and $(\bar p,\bar q-s)$ if $v=0$ and if $0<s\le p$;
\item of types $(r,0)$ and $(\bar p-r,\bar q)$ if $u>0$ and if  $u+2\le r\le\bar p$;
\item of types $(0,s)$ and $(\bar p,\bar q-s)$ if $v>0$  and if $v+2\le s\le\bar q$.\end{enumerate}
\item $(N,g_N)$ is not Jordan Osserman of type $(r,s)$ otherwise.\end{enumerate}\end{enumerate}
\end{theorem}

The metrics corresponding to the tensors $\psi_f$ are realizable as hypersurfaces
$S_f$ in a flat space and were studied in
\cite{refGIZ}. As the Hessian
$H(f)$ is the second fundamental form of $S_f$, the condition $H(f)>0$ of Theorem \ref{Lem1.5} (1) is simply the condition
that $S_f$ is concave. By using assertions (2) and (3) of Theorem \ref{Lem1.5}, we may perturb these
metrics to construct examples of higher order Jordan Osserman manifolds which are not hypersurfaces. Therefore, Lemma
\ref{Lem1.5} shows that there are manifolds described by Theorems \ref{Thm1.6} and \ref{Thm1.7} which are neither locally
homogeneous nor of constant sectional curvature for generic $\psi\in\Psi$.

\section{The Curvature Tensor of $(M,g_\psi)$}\label{Sect2}

If $Z_1$, $Z_2$, and $Z_3$ are coordinate vector fields, then we have:
\begin{equation}\label{eqn2.1}
(\nabla_{Z_1}Z_2,Z_3)=\textstyle\frac12\{Z_1(Z_2,Z_3)+Z_2(Z_1,Z_3)-Z_3(Z_1,Z_2)\}.
\end{equation}
We use equation (\ref{eqn2.1}) to see that $(\nabla_{Z_1}\partial_i^y,Z_3)=0$ for all $Z_1$ and $Z_3$. Since the
inner product is non-degenerate, this implies that
$\nabla\partial_i^y=0$. Consequently,
$$R(Z_1,Z_2)\partial_i^y=
(\nabla_{Z_1}\nabla_{Z_2}-\nabla_{Z_2}\nabla_{Z_1}-\nabla_{[Z_1,Z_2]})\partial_i^y=0.$$
Thus $R(Z_1,Z_2,Z_3,Z_4)=0$ if $Z_3\in\mathcal{Y}$. We use the curvature symmetries
\begin{eqnarray*}
&&\phantom{=-}R(Z_1,Z_2,Z_3,Z_4)=R(Z_3,Z_4,Z_1,Z_2)\\&&=-R(Z_4,Z_3,Z_1,Z_2)=-R(Z_1,Z_2,Z_4,Z_3)\end{eqnarray*}
 to see that
$R(Z_1,Z_2,Z_3,Z_4)=0$ if any of the $Z_i\in\mathcal{Y}$; this proves Lemma \ref{Lem1.3} (1).

By equation (\ref{eqn2.1}),
$(\nabla_{\partial_i^x}\partial_j^x,\partial_k^x)=\textstyle\frac12(\psi_{ik/j}+\psi_{jk/i}-\psi_{ij/k})$
and $(\nabla_{\partial_i^x}\partial_j^x,\partial_k^y)=0$.
Consequently,
$$\nabla_{\partial_i^x}\partial_j^x
   =\textstyle\frac12\sum_k(\psi_{ik/j}+\psi_{jk/i}-\psi_{ij/k})\partial_k^y.$$
We complete the proof of Lemma \ref{Lem1.3} by computing:
\begin{eqnarray*}
&&\ \ \ R(\partial_i^x,\partial_j^x)\partial_k^x=
   (\nabla_{\partial_i^x}\nabla_{\partial_j^x}-\nabla_{\partial_j^x}\nabla_{\partial_i^x})\partial_k^x\\
&&=\textstyle\frac12\nabla_{\partial_i^x}\{(\psi_{jl/k}+\psi_{kl/j}-\psi_{jk/l})\partial_l^y\}
   -\frac12\nabla_{\partial_j^x}\{(\psi_{il/k}+\psi_{kl/i}-\psi_{ik/l})\partial_l^y\}\\
&&=\textstyle\frac12\{\partial_i(\psi_{jl/k}+\psi_{kl/j}-\psi_{jk/l})-
  \partial_j(\psi_{il/k}+\psi_{kl/i}-\psi_{ik/l})\}\partial_l^y\\
&&=\textstyle\frac12(\psi_{jl/ki}+\psi_{kl/ji}-\psi_{jk/il}-\psi_{il/jk}-\psi_{kl/ij}+\psi_{ik/jl})\partial_l^y.\qquad\qedbox
\end{eqnarray*}

\medbreak We now establish the assertions of Lemma \ref{Lem1.5}. Let $f\in C^\infty(\mathbb{R}^p)$. We then set
$H_{ij}=\partial_i^x\partial_j^xf$ and
$\psi_{ij}:=\partial_i^xf\cdot\partial_j^xf$. By Lemma \ref{Lem1.3},
\begin{eqnarray}
(R(\partial_i^x,\partial_j^x)\partial_k^x,\partial_l^x)
&=&-\textstyle\frac12(\psi_{jk/il}+\psi_{il/jk}-\psi_{ik/jl}-\psi_{jl/ik})\nonumber\\
&=&H_{il}H_{jk}-H_{ik}H_{jl},\quad\text{so}\nonumber\\
(J(X_1)X_2,X_3)&=&H(X_1,X_1)H(X_2,X_3)-H(X_1,X_2)H(X_1,X_3).\label{eqn2.2}\end{eqnarray}
Suppose $H$ is positive definite. If $0\ne X_1\in\mathcal{X}$, let 
$$\mathcal{Z}=\{X_2\in\mathcal{X}:H(X_1,X_2)=0\}.$$ 
By equation
(\ref{eqn2.2}),
$(J(X_1)\cdot,\cdot)$ is positive definite on $\mathcal{Z}$ because
$$(J(X_1)X_2,X_3)=H(X_1,X_1)H(X_2,X_3)\quad\text{ for }\quad Z_2,Z_3\in\mathcal{Z}.$$ Since $\rank(J(X_1))\le p-1$, we have that
$J(X_1)$ is positive semi-definite of rank
$p-1$. This proves Lemma \ref{Lem1.5} (1).
Asssertion (2) of Lemma \ref{Lem1.5} follows from Lemma \ref{Lem1.3} as
$$J_{a\psi}=aJ_\psi.$$

 Let $\psi_1,\psi_2\in\Psi$, let $a_1>0$, and let $a_2>0$. By Lemma \ref{Lem1.3},
\begin{eqnarray*}
&&R_{a_1\psi_1+a_2\psi_2}(\partial_i^x,\partial_j^x)\partial_k^x=
a_1R_{\psi_1}(\partial_i^x,\partial_j^x)\partial_k^x+a_2R_{\psi_2}(\partial_i^x,\partial_j^x)\partial_k^x,\\
&&J_{a_1\psi_1+a_2\psi_2}(X_1)=a_1J_{\psi_1}(X_1)+a_2J_{\psi_2}(X_1).\end{eqnarray*}
Thus $J_{a_1\psi_1+a_2\psi_2}(X_1)$ is positive semi-definite of rank at least $p-1$ for $0\ne X_1\in\mathcal{X}$.  Since
$\rank(J_\psi(X_1))\le p-1$ for any $0\ne X_1\in\mathcal{X}$, $\rank(J_{a_1\psi_1+a_2\psi_2}(X_1))=p-1$ so
$a_1\psi_1+a_2\psi_2\in\Psi$. This establishes Lemma \ref{Lem1.5} (3).

 We complete the proof of Lemma \ref{Lem1.5} by studying small compactly supported perturbations. Previously, we have suppressed
the point $(x,y)$ of $M$ from the notation. It is now necessary to introduce it explicitly. Thus we will denote the Jacobi
operator by $J(x,y;Z)$ for $(x,y)\in M$ and
$Z\in\mathbb{R}^{2p}$. By Lemma \ref{Lem1.3}, $J(x,y;Z)=J(x;X)$ for $X=\pix Z$.

We must introduce an auxiliary positive definite quadratic form $\langle\cdot,\cdot\rangle$ on
$\mathbb{R}^p$. Let $S^{p-1}$ be the associated unit sphere. Let $SK:=K\times S^{p-1}$. Then
$$T(SK)=\{(x;X_2,X_3)\in(\mathbb{R}^p)^3:x\in K, \langle X_2,X_2\rangle=1,\langle X_2,X_3\rangle=0\}.$$
Let $\psi\in\Psi$. Then $J_\psi(x;X_2)$ defines a bilinear symmetric quadratic form $\mathcal{J}(x;X_2)$ on $T_{(x;X_2)}(SK)$;
$J_\psi(x;X_2)$ is positive semi-definite of rank $p-1$ if and only if $\mathcal{J}(x;X_2)$ is positive definite. Since $SK$ is
compact and since we are considering compact perturbations, Lemma \ref{Lem1.5} now follows from Lemma \ref{Lem1.3}.\hfill\qedbox

\section{The higher order Jacobi operator of $(M,g_\psi)$}\label{Sect4}

By Lemma \ref{Lem1.3}, $J(Z_1)Z_3=R(Z_3,Z_1)Z_1=0$ if $Z_1\in\mathcal{Y}$; thus $\mathcal{Y}\subset\ker(J(Z_1))$.
Furthermore, $\range(J(Z_2))\subset\Pspan\{R(\partial_i^x,\partial_j^x)\partial_k^x\}\subset\mathcal{Y}$. Thus
\begin{equation}J(Z_1)J(Z_2)=0.\label{eqn3.1}\end{equation}
Let $\{e_i\}$ be an orthonormal basis for $\pi\in\Gr_{r,s}(M,g_\psi)$. We then have
$$\textstyle J(\pi)^2=\sum_{i,j}(e_i,e_i)(e_j,e_j)J(e_i)J(e_j)=0.$$
Therefore, $\spec(J(\pi))=\{0\}$ so $(M,g_\psi)$ is
 Osserman of type $(r,s)$ and hence is $k$ Osserman for any admissible $k$. This proves Theorem \ref{Thm1.6} (1).

\medbreak Recall from equation (\ref{eqn1.1}) that $\pix$ is projection on $\mathcal{X}$. We begin the
proof of Theorem
\ref{Thm1.6} (2) with a technical Lemma:
\begin{lemma}\label{Lem3.1} Let $\psi\in\Psi$, let $\pi:=\Pspan\{Z_1,...,Z_r\}$ and let $J:=J(Z_1)+...+J(Z_r)$.
\begin{enumerate}
\item If $\dim(\pix\pi)=0$, then $\rank(J)=0$.
\item If $\dim(\pix\pi)=1$, then $\rank(J)=p-1$.
\item If $\dim(\pix\pi)>1$, then $J$ is positive semi-definite of rank $p$.
\end{enumerate}
\end{lemma}

\begin{proof} Let $X_i=\pix(Z_i)$; $J(Z_i)=J(X_i)$. Thus if $X_i=0$, then $J(Z_i)=0$. Suppose $X_i\ne0$. Since $\psi\in\Psi$, we
have that $\rank\{J(X_i)\}=p-1$. Because
$J(X_i)X_i=0$, 
\begin{equation}\ker(J(X_i))=\mathbb{R}\cdot X_i\oplus\mathcal{Y}\quad\text{ if }\quad X_i\ne0.\label{eqn3.2}\end{equation}
If $\dim(\pix\pi)=0$, then $X_i=0$ for all $i$ and $J=0$. Suppose $\dim(\pix \pi)=1$. We may suppose without loss of generality that
$X_1\ne0$ and let
$X_i=c_iX_1$ for
$i\ge1$. Then $J=(\sum_ic_i^2)J(X_1)$ has rank $p-1$. 

Suppose $\dim(\pix \pi)\ge2$. We suppose, without loss of generality, that $\{X_1,X_2\}$ is a linearly independent set.
By Lemma \ref{Lem1.3}, $\range(J)\subseteq\mathcal{Y}$. Therefore, we have that $\rank(J)\le\dim(\mathcal{Y})=p$. Conversely,
$J(X_i)$ is positive semidefinite as
$$(JZ,Z)=(J(X_1)Z,Z)+...+(J(X_r)Z,Z)\ge0.$$
Furthermore, equality holds if and only if $(J(X_i)\pix Z,\pix Z)=0$ for $1\le i\le
p$. By equation (\ref{eqn3.2}), this means that $\pix Z$ is a multiple of $X_i$ for $i=1,2$. Consequently, $\pix Z=0$ so
$Z\in\mathcal{Y}$. Therefore, $\rank(J)=p$.
\end{proof}

\medbreak Let $\psi\in\Psi$. Let $\{Z_1,...,Z_r\}$ be an orthonormal basis for a
spacelike subspace
$\pi\in\Gr_{0,r}(M,g)$. Let $X_i:=\pix Z_i$. Because $\pi$ is spacelike, $\pi\cap\mathcal{Y}=\{0\}$. Since $\ker\pix=\mathcal{Y}$,
$\dim(\pix\pi)=r$.  We have
$$J(\pi)=J(X_1)+...+J(X_r).$$
Since $J(\pi)^2=0$, the Jordan normal form of $J(\pi)$ is determined by the rank. We apply Lemma \ref{Lem3.1} to see that the rank
is
$p-1$ if $r=1$ and that the rank is $p$ if $r>1$. This shows $(M,g_\psi)$ is Jordan Osserman of type $(0,r)$. One shows similarly
that $(M,g_\psi)$ is Jordan Osserman of type $(r,0)$. Lemma \ref{Lem1.3} then shows that $(M,g_\psi)$ is Jordan Osserman of types
$(p,p-r)$ and $(p-r,p)$. Thus Theorem \ref{Thm1.6} (2a) is established.

\medbreak To prove Theorem \ref{Thm1.6} (2b), it is convenient to define the elements:
$$\tilde X_i:=\partial_i^x-\textstyle\frac12\sum_j\psi_{ij}\partial_j^y.$$
Note that $\pix\tilde X_i=\partial_i^x$ and that $\{\tilde X_1,...,\tilde X_p,\partial_1^y,...,\partial_p^y\}$ is a
{\it hyperbolic basis}, i.e.
\begin{equation}
  (\tilde X_i,\tilde X_j)=(\partial_i^y,\partial_j^y)=0\quad\text{and}\quad\ 
(\tilde X_i,\partial_j^y)=\delta_{ij}\quad\text{for}\quad1\le i,j\le p.\label{eqn3.3}\end{equation}

Let $0<r<p$ and $0<s<p$. We must show that $(M,g_\psi)$ is not Jordan Osserman of type $(r,s)$. By Theorem
\ref{Thm1.1} (2), we may assume $r+s\le p$. We assume $0<r\le s<p$ as the situation when $0<s\le r<p$ is similar.
We distinguish two cases:
\medbreak 1) Suppose $s=1$. Then $r=1$. We use equation (\ref{eqn3.3}) to define the following subspaces of type $(1,1)$
with the indicated orthonormal bases and Jacobi operators:
$$\begin{array}{ll}
\qquad\quad\pi_1:=
\textstyle\Pspan\{\tilde X_1-\frac12\partial_1^y,\tilde X_1+\frac12\partial_1^y\},&
J(\pi_1)=-J(\tilde X_1)+J(\tilde X_1),\\
\qquad\quad\pi_2:=\Pspan\{\varepsilon\tilde X_1-\frac12\varepsilon^{-1}\partial_1^y,
\tilde X_2+\frac12\partial^y_2\},&
J(\pi_2)=-\varepsilon J(\tilde X_1)+J(\tilde X_2).\end{array}$$
As $J(\pi_1)=0$,
$\rank(J(\pi_1))=0$. 
Since $\rank(J(\tilde X_2))=\rank(J(\partial_2^x))=p-1$,
$\rank(J(\pi_2))\ge\rank(J(\tilde
X_2))\ge p-1$ if $\varepsilon$ is small. Consequently, $(M,g_\psi)$ is not Jordan Osserman of type $(1,1)$.

\medbreak 2) Suppose $1\le r\le s< p$, $s\ge2$, and $r+s\le p$; necessarily $p\ge3$. For
$\alpha\ne0$ and $\beta\ne0$, we use equation (\ref{eqn3.3}) to define timelike subspaces
$\pi^-(\alpha)$ and spacelike subspaces $\pi^+(\beta)$ with the indicated orthonormal bases and Jacobi operators:
\begin{eqnarray*}
&&\textstyle\pi^-(\alpha):=\Pspan\{\alpha\tilde X_1-\frac12\alpha^{-1}\partial_1^y,
...,\alpha\tilde X_r-\frac12\alpha^{-1}\partial^y_r\},\\
&&J^-_\alpha:=-\alpha\{J(\partial_1^x)+...+J(\partial_r^x)\},\\
&&\pi^+(\beta):=\textstyle
\Pspan\{\beta\tilde X_{r+1}+\frac12\beta^{-1}\partial^y_{r+1},...,
\beta\tilde X_{r+s}+\frac12\beta^{-1}\partial^y_{r+s}\},\\
&&\textstyle J_\beta^+:=\beta\{J(\partial_{r+1}^x)+...+J(\partial_{r+s}^x)\}.\end{eqnarray*}
Let $\pi(\alpha,\beta):=\pi^-(\alpha)\oplus\pi^+(\beta)\in Gr_{r,s}(M,g)$. The associated Jacobi
operator $J(\pi(\alpha,\beta))=J_\alpha^-+J_\beta^+$.
Since $\range(J)\subset\mathcal{Y}$, one sees that
$$\rank(J(\pi(\alpha,\beta)))\le p\quad\text{ for all }\quad\alpha,\beta.$$
We have $\pix(\tilde X_i)=\partial_i^x$. We apply Lemma \ref{Lem3.1}; $J(\cdot)$ is supported on $\mathcal{X}$. Since
$s\ge2$, $J^+_\beta$ is positive semi-definite of rank
$p$. Thus $\rank\{J(1,\beta)\}=p$ for $\beta$ large.
Note that $J^-_\alpha$
is negative semi-definite of rank at least $p-1$. Thus for
$\alpha$ large,
$J(\alpha,1)$ determines a quadratic form of signature $(u,v)$ for $u\ge p-1$ and $u+v\le p$. Thus by continuity,
there must exist
$(\alpha,\beta)$ with $\alpha\ne0$ and $\beta\ne0$ so $J(\pi(\alpha,\beta))$ determines a degenerate quadratic form on
$\mathcal{X}$. For such values of $(\alpha,\beta)$,
$\rank\{J(\pi(\alpha,\beta))\}<p$. This shows that $\rank\{J(\pi(\alpha)\oplus\pi(\beta))\}$ is not constant and hence
$(M,g_\psi)$ is not Jordan Osserman of type $(r,s)$. The proof of Theorem
\ref{Thm1.6} is now complete.\hfill\qedbox

\section{The higher order Jacobi operator of $(N,g_N)$}\label{Sect5}

Let
$N:=M\times\mathbb{R}^{(u,v)}$ with the product metric. Decompose 
$TN=\mathcal{X}\oplus\mathcal{Y}\oplus T\mathbb{R}^{(u,v)}$ and
let $\pix :TN\rightarrow\mathcal{X}$ be the associated projection. Then
$$J(U_1)U_2=J(\pix U_1)\pix U_2\quad\text{ for any }\quad U_1,U_2\in TN.$$
Let $\{U_\nu\}$ be an orthonormal basis for a non-degenerate subspace $\pi\in Gr_{r,s}(N)$. Let $\varepsilon_\nu:=(U_\nu,U_\nu)=\pm1$.
Then $J(\pi)=\sum_\nu\varepsilon_\nu J(\pix(U_\nu))$ so equation (\ref{eqn3.1}) yields:
$$\textstyle
J(\pi)^2=\sum_{\mu,\nu}\varepsilon_\nu\varepsilon_\mu J(\pix(U_\nu))J(\pix(U_\mu))=0.$$
This shows that $(N,g_N)$ is Osserman of type $(r,s)$ for every admissible $(r,s)$. Therefore $(N,g_N)$ is $k$ Osserman for every
admissible $k$ which completes the proof of Theorem
\ref{Thm1.7} (1). 

\medbreak Since $J(\pi)^2=0$,
$(N,g_N)$ will be Jordan Osserman of type $(r,s)$ if and only if $\rank(J(\pi))$ is constant on $\Gr_{r,s}(N,g_N)$. 
We use this observation to prove Theorem \ref{Thm1.7} (2). Let $\{U_1,...,U_r\}$ be an orthonormal basis for a timelike subspace
$\pi$ in $\Gr_{r,0}(N,g_N)$. We have
\begin{eqnarray*}
&&\pix(\pi)=\Pspan\{\pix(U_1),...,\pix(U_r)\},\text{ and}\\
&&J(\pi)=-J(\pix(U_1))-...-J(\pix(U_r)).\end{eqnarray*}
Lemma \ref{Lem3.1} implies:
\begin{equation}\rank(J(\pi))=\left\{\begin{array}{ll}
0&\text{if }\dim\pix(\pi)=0,\\
p-1&\text{if }\dim\pix(\pi)=1,\\
p&\text{if }\dim\pix(\pi)\ge2.\end{array}\right.\label{eqn4.1}\end{equation}

 Suppose that $u=0$. Then $\dim(\pix(\pi))=r$ is independent of $\pi$ and, by equation (\ref{eqn4.1}),
$(N,g_N)$ is Jordan Osserman of type $(r,0)$. Dually, by Theorem \ref{Thm1.1}, $(N,g_N)$ is Jordan Osserman of type $(p-r,p+v)$.
This proves Theorem \ref{Thm1.7} (2a-i).

Suppose $r\ge u+2$ and $u>0$. Then $\dim\pix(\pi)\ge r-u\ge2$ and hence
$\rank(J(\pi))=p$. Thus $(N,g_N)$ is Jordan Osserman of type $(r,0)$ and dually, it is Jordan Osserman of type $(p-r,v+p)$. This
completes the proof of assertion (2a-iii).
Assertions (2a-ii) and (2a-iv) are proved similarly.

\medbreak The proof of assertion (3) decomposes into several cases. Suppose $u>0$ and $r\le u+1$. We must show that $(N,g_N)$ is
not Jordan Osserman of type
$(r,0)$; we may then use duality to see that $(N,g_N)$ is not Jordan Osserman of type $(\bar p-r,\bar q)$. Let
$\{V_1^-,...,V_u^-\}\subset TN$ and $\{Z_1^-,...,Z_p^-\}\subset TM$ be orthonormal sets of timelike vectors. Let
$$\pi(a,b):=\Pspan\{Z_1^-,...,Z_a^-,V_1^-,...,V_b^-\}.$$
 Suppose first that $0<r\le u$. Then
\begin{eqnarray*}
&&\dim\pix(\pi(1,r-1))=1\quad\text{ so }\quad\rank(J(\pi(1,r-1))=p-1,\\
&&\dim\pix(\pi(0,r))=0\phantom{......}\quad\text{ so }\quad\rank(J(\pi(0,r))=0.\end{eqnarray*}
Therefore $(N,g_N)$ is not Jordan Osserman of type $(0,r)$. Suppose next $r=u+1$. We then have $r\ge2$ and
\begin{eqnarray*}
&&\dim\pix(\pi(2,u-1))=2\quad\text{ so }\quad\rank(J(\pi(2,u-1)))=p,\\
&&\dim\pix(\pi(1,u))=1\quad\phantom{......}\text{ so }\quad\rank(J(\pi(1,u)))=p-1.\end{eqnarray*}
Thus $(N,g_N)$ is not Jordan Osserman of type $(0,r)$.
The case $v>0$ and $r<v+2$ is similar and is omitted in the interests of brevity.

Finally, suppose $1\le r\le \bar p-1$ and $1\le s\le \bar q-1$. Let 
$$\{V_1^-,...,V_u^-,V_1^+,...,V_v^+\}$$ be an orthonormal basis for
$\mathbb{R}^{(u,v)}$. We define maps
$$T_{a,b}\pi:=\pi\oplus\Pspan\{V_1^-,...,V_a^-,V_1^+,...,V_b^+)$$
from $\Gr_{\alpha,\beta}(M,g_\psi)$ to $\Gr_{\alpha+a,\beta+b}(N,g_N)$.
We then have 
\begin{equation}J_N(T_{a,b}\pi)=J_M(\pi)\quad\text{ for all }\quad\pi\in\Gr_{\alpha,\beta}(M,g_\psi).\label{eqn4.2}\end{equation}
Suppose $(N,g_N)$ is Jordan Osserman of type $(r,s)$. Expand 
\begin{eqnarray*}
&&r=\alpha+a,\quad\text{ where }\quad1\le\alpha\le p-1\quad\text{ and }\quad0\le a\le u;\\
&&s=\beta+b,\phantom{.}\quad\text{ where }\quad1\le\beta\le p-1\quad\text{ and }\quad0\le b\le v.\end{eqnarray*}
If $(N,g_N)$ is Jordan Osserman of type $(r,s)$, then we
may use equation (\ref{eqn4.2}) to see that $(M,g_\psi)$ is Jordan Osserman of type $(\alpha,\beta)$. This contradicts Theorem
\ref{Thm1.6} and thereby completes the proof of Theorem \ref{Thm1.7}.
\hfill\qedbox

\bigbreak\noindent{\bf Acknowledgments:} 
Stavrov \cite{refSt} constructed algebraic curvature tensors which are Jordan Osserman for some, but not all, values of
$(r,s)$. Our analysis in the geometric context is motivated at least in part by her analysis in the algebraic setting and it is a
pleasure to acknowledge helpful conversations with her concerning these matters.

The research of P. Gilkey was partially supported by the NSF (USA) and
the MPI (Leipzig); the research of R. Ivanova and T. Zhang was partially supported by the NSF (USA).


\begin{thebibliography}{AAA}

\bibitem{refBokanBlazicGilkey} N. Bla\v zic, N. Bokan, and P. Gilkey, A
    note on Osserman Lorentzian manifolds, {\it Bulletin of the London Math Society},  {\bf 29},
       (1997), 227--230.

\bibitem{refBBGZ} N. Bla\v zi\'c, N. Bokan, P. Gilkey and Z. Raki\'c,
     Pseudo-Riemannian Osserman manifolds, {\it J. Balkan Soc.
     of Geometers}, {\bf l2}, (1997), 1--12.

\bibitem{refBCG} A. Bonome, P. Castro, E. Garcia-Rio, 
Four-Dimensional Generalized Osserman Manifolds,
{\it Classical and Quantum Gravity}, {\bf 18} (2001), 4813--4822.



\bibitem{refChia} Q.-S. Chi, A curvature characterization of certain locally
rank-one symmetric spaces, {\it J. Differ. Geom.}, {\bf 28}, (1988), 187--202.

\bibitem{refGra} E. Garc\'ia-Rio, D. Kupeli, and M. V\' azquez-Abal,
On a problem of Osserman in Lorentzian geometry, {\it Diff. Geom. Appl.}, {\bf 7}, (1997),
85--100.

\bibitem{refGVV} E. Garc\'ia-Ri\'o, M. E. V\' azquez-Abal and
     R. V\' azquez-Lorenzo, Nonsymmetric Osserman pseudo-Riemannian manifolds,
     {\it Proc. Amer. Math. Soc.}, {\bf 126}, (1998), 2771--2778.

\bibitem{refGRKVL} E. Garci\'a-Ri\'o, D. Kupeli, and R. V\'azquez-Lorenzo, {\bf Osserman Manifolds in Semi-Riemannian
Geometry}, Lecture notes in Mathematics, Springer Verlag, (2002), ISBN 3-540-43144-6.

\bibitem{refGiA} P. Gilkey,  Algebraic curvature tensors which are $p$ Osserman, {\it Diff.
Geometry and Appl.}, {\bf 14}, (2001), 297--311.

\bibitem{refGil} ---, {\bf Geometric Properties of Natural Operators Defined by the Riemann
Curvature Tensor}, World Scientific Press (2001), ISBN 981-02-04752-4.

\bibitem{refGiIv} P. Gilkey and R. Ivanova, The Jordan normal form of higher order Osserman algebraic curvature
tensors,  {\it Comment. Math. Univ. Carolinae}, {\bf 43}, (2002) 231--242.

\bibitem{refGIZ} P. Gilkey, R. Ivanova, and T. Zhang, Szab\'o Osserman IP Pseudo-Riemannian manifolds, preprint:
http://arXiv.org/abs/math.DG/0205085.

\bibitem{refGSV} P. Gilkey, G. Stanilov and V. Videv,  Pseudo-Riemannian
      manifolds whose generalized Jacobi operator
     has constant characteristic polynomial, {\it J. Geom.}, {\bf 62}, (1998), 144--153.

\bibitem{refGiSt} P. Gilkey and I. Stavrov,  Curvature tensors
   whose Jacobi or Szab\'o operator is nilpotent
    on null vectors, {\it Bull. London Math. Soc.}, to appear.

\bibitem{refNik} Y. Nikolayevsky, Osserman Conjecture in dimension $n \ne 16$, preprint:\newline
http://arXiv.org/abs/math.DG/0204258.


\bibitem{refOss} R. Osserman, Curvature in the eighties, {\it Amer. Math.
    Monthly}, {\bf97}, (1990), 731--756.

\bibitem{refSV} G. Stanilov and V. Videv, Four dimensional
    pointwise Osserman manifolds, {\it Abh. Math. Sem. Univ. Hamburg}, {\bf 68}, (1998),
    1--6.

\bibitem{refSt} I. Stavrov, Ph. D. Thesis, University of Oregon (2003).

\end{thebibliography}
\end{document}